\documentclass{amsart}
\usepackage{latexsym}
\usepackage{amssymb}
\usepackage{amsmath}

\begin{document}

%%%%%%%%%%%%%%%%%%%%%%Definitions%%%%%%%%%%%%%%%%%%%%%%%%%%%%%%%%%%%%%%%%%%%

\newtheorem{theorem}{Theorem}
\newtheorem{problem}{Problem}
\newtheorem{definition}{Definition}
\newtheorem{lemma}{Lemma}
\newtheorem{proposition}{Proposition}
\newtheorem{corollary}{Corollary}
\newtheorem{example}{Example}
\newtheorem{conjecture}{Conjecture}
\newtheorem{algorithm}{Algorithm}
\newtheorem{exercise}{Exercise}
\newtheorem{remarkk}{Remark}

\newcommand{\be}{\begin{equation}}
\newcommand{\ee}{\end{equation}}
\newcommand{\bea}{\begin{eqnarray}}
\newcommand{\eea}{\end{eqnarray}}
\newcommand{\beq}[1]{\begin{equation}\label{#1}}
\newcommand{\eeq}{\end{equation}}
\newcommand{\beqn}[1]{\begin{eqnarray}\label{#1}}
\newcommand{\eeqn}{\end{eqnarray}}
\newcommand{\beaa}{\begin{eqnarray*}}
\newcommand{\eeaa}{\end{eqnarray*}}
\newcommand{\req}[1]{(\ref{#1})}

\newcommand{\lip}{\langle}
\newcommand{\rip}{\rangle}

\newcommand{\uu}{\underline}
\newcommand{\oo}{\overline}
\newcommand{\La}{\Lambda}
\newcommand{\la}{\lambda}
\newcommand{\eps}{\varepsilon}
\newcommand{\om}{\omega}
\newcommand{\Om}{\Omega}
\newcommand{\ga}{\gamma}
\newcommand{\rrr}{{\Bigr)}}
\newcommand{\qqq}{{\Bigl\|}}

\newcommand{\dint}{\displaystyle\int}
\newcommand{\dsum}{\displaystyle\sum}
\newcommand{\dfr}{\displaystyle\frac}
\newcommand{\bige}{\mbox{\Large\it e}}
\newcommand{\integers}{{\Bbb Z}}
\newcommand{\rationals}{{\Bbb Q}}
\newcommand{\reals}{{\rm I\!R}}
\newcommand{\realsd}{\reals^d}
\newcommand{\realsn}{\reals^n}
\newcommand{\NN}{{\rm I\!N}}
\newcommand{\degree}{{\scriptscriptstyle \circ }}
\newcommand{\dfn}{\stackrel{\triangle}{=}}
\def\complex{\mathop{\raise .45ex\hbox{${\bf\scriptstyle{|}}$}
     \kern -0.40em {\rm \textstyle{C}}}\nolimits}
\def\hilbert{\mathop{\raise .21ex\hbox{$\bigcirc$}}\kern -1.005em {\rm\textstyle{H}}} %Hilbert space
\newcommand{\RAISE}{{\:\raisebox{.6ex}{$\scriptstyle{>}$}\raisebox{-.3ex}
           {$\scriptstyle{\!\!\!\!\!<}\:$}}} % >< one above each other

\newcommand{\hh}{{\:\raisebox{1.8ex}{$\scriptstyle{\degree}$}\raisebox{.0ex}
           {$\textstyle{\!\!\!\! H}$}}}

\newcommand{\OO}{\won}
\newcommand{\calA}{{\mathcal A}}
\newcommand{\calB}{{\mathcal B}}
\newcommand{\calC}{{\cal C}}
\newcommand{\calD}{{\cal D}}
\newcommand{\calE}{{\cal E}}
\newcommand{\calF}{{\mathcal F}}
\newcommand{\calG}{{\cal G}}
\newcommand{\calH}{{\cal H}}
\newcommand{\calK}{{\cal K}}
\newcommand{\calL}{{\mathcal L}}
\newcommand{\calM}{{\cal M}}
\newcommand{\calO}{{\cal O}}
\newcommand{\calP}{{\cal P}}
\newcommand{\calU}{{\mathcal U}}
\newcommand{\calX}{{\cal X}}
\newcommand{\calXX}{{\cal X\mbox{\raisebox{.3ex}{$\!\!\!\!\!-$}}}}
\newcommand{\calXXX}{{\cal X\!\!\!\!\!-}}
\newcommand{\gi}{{\raisebox{.0ex}{$\scriptscriptstyle{\cal X}$}
\raisebox{.1ex} {$\scriptstyle{\!\!\!\!-}\:$}}}
\newcommand{\intsim}{\int_0^1\!\!\!\!\!\!\!\!\!\sim}
\newcommand{\intsimt}{\int_0^t\!\!\!\!\!\!\!\!\!\sim}
\newcommand{\pp}{{\partial}}
\newcommand{\al}{{\alpha}}
\newcommand{\sB}{{\cal B}}
\newcommand{\sL}{{\cal L}}
\newcommand{\sF}{{\cal F}}
\newcommand{\sE}{{\cal E}}
\newcommand{\sX}{{\cal X}}
\newcommand{\R}{{\rm I\!R}}
\renewcommand{\L}{{\rm I\!L}}
\newcommand{\vp}{\varphi}
\newcommand{\N}{{\rm I\!N}}
\def\ooo{\lip}
\def\ccc{\rip}
\newcommand{\ot}{\hat\otimes}
\newcommand{\rP}{{\Bbb P}}
\newcommand{\bfcdot}{{\mbox{\boldmath$\cdot$}}}

\renewcommand{\varrho}{{\ell}}
\newcommand{\dett}{{\textstyle{\det_2}}}
\newcommand{\sign}{{\mbox{\rm sign}}}
\newcommand{\TE}{{\rm TE}}
\newcommand{\TA}{{\rm TA}}
\newcommand{\E}{{\rm E\,}}
\newcommand{\won}{{\mbox{\bf 1}}}
\newcommand{\Lebn}{{\rm Leb}_n}
\newcommand{\Prob}{{\rm Prob\,}}
\newcommand{\sinc}{{\rm sinc\,}}
\newcommand{\ctg}{{\rm ctg\,}}
\newcommand{\loc}{{\rm loc}}
\newcommand{\trace}{{\,\,\rm trace\,\,}}
\newcommand{\Dom}{{\rm Dom}}
\newcommand{\ifff}{\mbox{\ if and only if\ }}
\newcommand{\nproof}{\noindent {\bf Proof:\ }}
\newcommand{\remark}{\noindent {\bf Remark:\ }}
\newcommand{\remarks}{\noindent {\bf Remarks:\ }}
\newcommand{\note}{\noindent {\bf Note:\ }}

\newcommand{\boldx}{{\bf x}}
\newcommand{\boldX}{{\bf X}}
\newcommand{\boldy}{{\bf y}}
\newcommand{\boldR}{{\bf R}}
\newcommand{\uux}{\uu{x}}
\newcommand{\uuY}{\uu{Y}}

\newcommand{\limn}{\lim_{n \rightarrow \infty}}
\newcommand{\limN}{\lim_{N \rightarrow \infty}}
\newcommand{\limr}{\lim_{r \rightarrow \infty}}
\newcommand{\limd}{\lim_{\delta \rightarrow \infty}}
\newcommand{\limM}{\lim_{M \rightarrow \infty}}
\newcommand{\limsupn}{\limsup_{n \rightarrow \infty}}

\newcommand{\ra}{ \rightarrow }

\newcommand{\ARROW}[1]
  {\begin{array}[t]{c}  \longrightarrow \\[-0.2cm] \textstyle{#1} \end{array} }

\newcommand{\AR}
 {\begin{array}[t]{c}
  \longrightarrow \\[-0.3cm]
  \scriptstyle {n\rightarrow \infty}
  \end{array}}

\newcommand{\pile}[2]
  {\left( \begin{array}{c}  {#1}\\[-0.2cm] {#2} \end{array} \right) }

\newcommand{\floor}[1]{\left\lfloor #1 \right\rfloor}

%for doing boldface subscripts etc., e.g. $G_{\mmbox{\boldx}}$
\newcommand{\mmbox}[1]{\mbox{\scriptsize{#1}}}

%fraction with round brackets
\newcommand{\ffrac}[2]
  {\left( \frac{#1}{#2} \right)}

\newcommand{\one}{\frac{1}{n}\:}
\newcommand{\half}{\frac{1}{2}\:}

\def\le{\leq}
\def\ge{\geq}
\def\lt{<}
\def\gt{>}

%qed
\def\squarebox#1{\hbox to #1{\hfill\vbox to #1{\vfill}}}
\newcommand{\nqed}{\hspace*{\fill}
           \vbox{\hrule\hbox{\vrule\squarebox{.667em}\vrule}\hrule}\bigskip}

\newcommand{\no}{\noindent}
\newcommand{\EE}{\mathbb{E}}
\newcommand{\RR}{\mathbb{R}}
\newcommand{\DD}{\mathbb{D}}
\newcommand{\F}{\mathcal{F}}
\newcommand{\W}{\mathbb{W}}
%%%%%%%%%%%%%%%%%%%%%%%%%%%%%%%%%%%%%%%%%%%%%%%%%%%%%%%%%%%%%%%%%%%%%%%%%%%%

\title{A CRITERIA OF STRONG H-DIFFERENTIABILITY}

\author{ K\'evin HARTMANN}
%\date{ }
\maketitle
\noindent
{\bf Abstract:}{\small{ We
give a criteria for a Malliavin differentiable function to be strongly H-differentiable.
}}\\

\vspace{0.5cm}

\noindent
Keywords: Wiener space, H-$C^1$, strong differentiability, Malliavin derivative\\

%\tableofcontents

\section{\bf{Introduction}}
\no Let $\W$ be the classical Wiener space and H
the associated Cameron-Martin space. A theory of a weak derivative over Wiener
functional with respect to H directions has long been developed (see
\cite{cours}, \cite{wat1}). More recently, \"Ust\"unel
and Zakai, in \cite{BOOK}, or Kusuoka, in \cite{ku1}   have studied
a
strong derivative for Wiener functional, using the Fr\'echet
differentiability on H. A Wiener functional f is H-continuous, or H-C, if
$h\mapsto f(w+h)$ is a.s. continuous on H, H-$C^1$ if $h\mapsto f(w+h)$ is
a.s. Fr\'echet differentiable on H with H-continuous derivative. H-$C^1$ function are very useful in the study of invertibility of perturbations of the identity of Wiener space. Indeed if u is a measurable H-$C^1$ function from $\W$ to H, $I_\W+u$ is invertible. This has been used by \"Ust\"unel to establish the following variational representation, where B is a Brownian motion
$$-\log\EE\left[e^{-f\circ B}\right]=\inf_u\EE\left[f\circ (B+u)+\frac{1}{2}\int_0^1|\dot{u}(s)|^2ds\right]$$
\no for some unbounded functions f.\newline
\no Of course it is way more difficult to establish that a function is H-$C^1$ than it is to establish it is weakly H-differentiable. In this paper, we give a criteria for a weakly H-differentiable function to be H-$C^1$, namely the weak H-derivative has to be a.s. uniformly continuous on every zero-centered ball of H.\newline\no First we recall the formal setting of weak and strong H-derivative, then we establish the criteria. Finally, we expand the criteria to higher order derivatives.

\section{\bf{Framework}}
\no
Set $n\in\NN$ and let $\W$ be the canonical Wiener space $C([0,1],\RR^n)$. Let H be the associated Cameron-Martin space
$$H=\left\{\int_0^.\dot{h}(s)ds,\dot{h}\in L^2([0,1])\right\}$$
\no and for $m\in\NN^*$, $B_m=\{h\in H, |h|_H\leq m\}$.
\no Denote $\mu$ the Wiener measure and W the coordinate process. W is a Brownian motion under $\mu$ and we denote $(\F_t)$ the canonical filtration of W completed with respect to $\mu$. Set Cyl the set of cylindrical functions
$$Cyl=\left\{F(W_{t_1},...,W_{t_p}),p\in\NN^*,F\in\mathcal{S}(\RR^n),0\leq t_1<...<t_p\leq 1\right\}$$
\no where $\mathcal{S}(\RR^n)$ denotes the set of Schwartz functions on $\RR^n$.\newline
\no For $f\in Cyl,w\in W$ and $h\in H$, we define
$$\nabla_h f(w)=\left.\frac{d}{d\lambda}f(w+\lambda h)\right|_{\lambda=0} $$
\no Riesz theorem enables us to consider $\nabla f$ as an element of H. For $1<p<\infty$, we define
$$|.|_{p,1}:f\in Cyl\mapsto|f|_{L^p(\mu)}+|\nabla f|_{L^p(\nu,H)}$$
\no $\nabla f$ is a closable operator and we define $\mathbb{D}_{p,1}$ the closure of Cyl for $|.|_{p,1}$.\newline
\no Let $\delta$ be the adjoint operator of $\nabla$ and $L_a^0(\mu,H)$ be the set of the element of $L^0(\mu,H)$ whose density are adapted to $(\F_t)$. $L^0_a(\mu,H)$ is a subset of the domain of $\delta$, and for any $u\in L^0_a(\mu,H)$
$$\delta u=\int_0^1\dot{u}(s)dW(s)$$
\no From now on for $u\in L^0_a(\mu,H)$, we will denote
$$\rho(\delta u)=\exp\left(\delta u-\frac{1}{2}\int_0^1|\dot{u}(s)|^2ds\right)$$
\no Now set X a separable Hilbert space and $(e_i)_{i\in\NN}$ an Hilbert base of X, define
$$Cyl(X)=\left\{\sum_{k=1}^p f_i e_{i_k}, p\in\NN^*,(i_k)\in\NN^p, (f_i)\in Cyl^p\right\}$$
\no If $f=\sum_{k=1}^p f_i e_{i_k}\in Cyl(X)$, we define
$$\nabla f (w)[h]=\sum_{k=1}^p \nabla_h f_i e_{i_k}$$
\no and $\nabla f$ is an element of $X\otimes H$.\newline
\no We define $|.|_{p,1}$, similarly as before and $\nabla$ is once again a closable operator, we define $\mathbb{D}_{p,1}(X)$ the closure of $Cyl(X)$ for $|.|_{p,1}$. This enables us to define $\nabla^p$ for $p\geq 1$ by recurrence, we denote
$$|.|_{p,k}:f\mapsto |f|_{L^p(\mu,X)}+\sum_{k=1}^p\left|\nabla^k f\right|_{L^p(\mu,X\otimes H^{\otimes^k})}$$
\no and $\mathbb{D}_{p,k}(X)$ the completion of $Cyl(X)$ for $|.|_{p,k}$.
\newline
\no Finally, we define the Ornstein-Uhlebeck semigroup $(P_t)$ as follow:  set $t>0$ and $f\in L^p(\mu,X)$ for some $p\geq 1$
$$P_tf:w\in\W\mapsto \int_\W f\left(e^{-t}w+\sqrt{1-e^{-2t}}y\right) \mu(dy)$$
\no We will need the following technical results concerning $P_t$:

\begin{proposition}
Set $t>0$ and $f\in L^1(\mu,X)$. For $h\in H$, we have $\mu$-a.s.
$$P_tf(w+h)=P_t\left((f\left((.+e^{-t}h\right)\right)(w)$$
\no If f belongs to some $\mathbb{D}_{p,1}(X)$, we have $\mu$-a.s.
$$P_t\nabla f=e^t\nabla P_t f$$
\end{proposition}

\nproof For the sake of simplicity we address the case $X=\RR$. The first assertion is an easy calculation.\newline
\no For the second one, set $h\in H$, we have
$$\nabla\rho(\delta h)=h\rho(\delta h)$$
\no and $$P_t\rho(\delta h)=\rho(\delta(e^{-t}h))$$
\no so
\beaa P_t\nabla\rho(\delta h)&=&h\rho(\delta(e^{-t} h))\\
&=&e^te^{-t}h\rho(\delta(e^{-t} h))\\
&=&e^t\nabla\rho(\delta(e^{-t} h))\\
&=&e^t\nabla P_t\rho(\delta h)\eeaa
\no and we conclude with density of the vector space generated by $\{\rho(\delta h),h\in H\}$ in $\mathbb{D}_{p,1}$.\nqed

\no For more details on this setting see \cite{cours} or \cite{wat1}.\newline
\no Now we give the definitions of strongly H-differentiable functions.

\begin{definition}
Set $u:\W\rightarrow X$ a measurable function. We say that\newline
\no (i) u is H-continuous (or H-C) if the map $h\mapsto u(w+h)$ is $\mu$-a.s. continuous on H.\newline
\no (ii) u is H-$C^1$ if the map $h\mapsto u(w+h)$ is $\mu$-a.s. Fr\'echet-differentiable and its Fr\'echet derivative $\nabla f$ is an H-continuous map from $\W$ to $X\otimes H$.\newline
\no (iii) Set $p\in\NN$, by recurrence, u is H-$C^p$ if it is H-$C^{p-1}$ and its derivative of order p-1 is an H-$C^1$ map from $\W$ to $X\otimes H^{\otimes^{p-1}}$
\end{definition}

\no We will need the following results concerning strong H-regularity for our main theorem, see \cite{BOOK} for their proof:

\begin{proposition}
Set $u:\W\mapsto X$ such that $\nabla^k u$ is well-defined for every $k\in\NN^*$. Assume there exists $p\in\NN^*$ such that for every $\lambda\in\RR_+$
$$\sum_{k=0}^\infty\frac{\lambda^k}{k!}\left|\nabla^k u\right|_{L^p(\mu,X\otimes H^{\otimes^k})}<\infty$$
\no Then $\mu$-a.s. for every $h\in H$
$$u(w+h)=\sum_{k=0}^\infty\frac{1}{k!}\nabla^k u(w)\left[h^{\otimes^k}\right]$$
\end{proposition}

\begin{proposition}
Set $f\in L^p(\mu,X)$ for some $p>1$, for every $t>0$ and $\lambda\in\RR_+$, we have
$$\sum_{k=0}^\infty\frac{\lambda^k}{k!}\left|\nabla^k P_t f\right|_{L^p(\mu,X\otimes H^{\otimes^k})}<\infty$$
\end{proposition}

\section{\bf{Main theorem}}

\begin{theorem}
\label{th}
\no Assume that $f:\W\rightarrow X$ is in $\mathbb{D}_{p,1}(X)$ for
some $p>1$. Assume that
$h\mapsto \nabla f(w+h)$ uniformly continuous on every n-ball of H.
\no Then f is $H-C^1$ and its H-derivative is $\nabla f$.
\end{theorem}

\nproof
\no The hypothesis implies that $h\mapsto \nabla f(w+h)$ is
separable so the uniform continuity hypothesis can be written:
$$\lim_{\epsilon\rightarrow 0}\sup_{h,k\in B_n, |h-k|_H\leq\epsilon}
|\nabla f(w+h)-\nabla f(w+k)|_{X\otimes H}=0\;\;a.s.$$
As we just stated we can set $A\subset \W$ of full measure such
that for every $w\in \W$ $h\mapsto\nabla f(w+h)$ is continuous.\newline
\no Set $s>0$ and $h\in H$. We know the action of $P_s$ over the weak derivative:
\beaa P_s\nabla f(w+h)=e^s\nabla P_s f(w+h)\;\;a.s.\eeaa
\no We also have:
\beaa P_s\nabla f(w+h)=P_s(\nabla f(.+e^{-s}h))(w)\;\;a.s.\eeaa
\no Since both terms are analytic, the set on which these equalities
hold does not depend on h. Now we denote, for $m,n\in\NN^*$:
$$\theta_{nm}(w)=\sup_{h,k\in B_n,|h_k|_H\leq\frac{1}{m}} \left|\nabla
  f(w+h)-\nabla f(w+k)\right|_{X\otimes H}$$
\no Observe that for $h,k\in B_n$ verifying $|h-k|_H<\frac{1}{m}$, we
have:
$$\left|P_s(\nabla f(.+e^{-s}h))(w)-P_s(\nabla f(.+e^{-s}k))(w)\right|_{X\otimes
  H}\leq P_s \theta_{nm}(w)\;\;a.s.$$
Since both terms have analytic modifications, the set of w on which this
inequality stands is independent of h and k.\newline
\no Set $(s_i)$ a sequence decreasing towards 0 and $H_0$ a countable
dense subset of H. We define:

\beaa A'&=& A\\
&\cap&\left\{w\in \W: P_{s_i}\nabla f(w+h)=e^{s_i}\nabla
  P_{s_i}f(w+h),\forall h\in H,\forall i\in\NN\right\}\\
&\cap&\{w\in \W: P_{s_i}\nabla f(w+h)=P_{s_i}(\nabla
f(.+e^{-s_i}h))(w), \forall h\in H, \forall i\in\NN\}\\
&\cap& \left\{w\in \W:\lim_{i\rightarrow\infty}
  P_{s_i}\nabla f(w+h)=\nabla f(w+h)\;\;\forall h\in H_0\right\}\\
&\cap& \left\{w\in \W:\lim_{i\rightarrow\infty}
  P_{s_i}\theta_{nm}(w)=\theta_{nm}(w),\forall n,m\in\NN\right\}\\
&\cap& \left\{w\in \W: \left|P_{s_i}(\nabla
    f(.+e^{-s_i}h))(w)-P_{s_i}(\nabla
    f(.+e^{-s_i}k))(w)\right|_{X\otimes H}\leq P_{s_i}\theta_{nm}(w),\right.\\
&&\left.\forall h,k\in B_n, |h-k|_H<1/m,\forall i\in\NN \right\}\\
&\cap& \left\{w\in \W: \lim_{m\rightarrow\infty}\theta_{nm}(w)=0,
  \forall n\in\NN \right\}\\
&\cap& \left\{w\in \W:P_{s_i}\nabla f(w+h)=P_{s_i}\nabla
  f(w)+\sum_{k=1}^\infty\frac{1}{k!}\nabla^k
  P_{s_i}\nabla f(w)\left[h^{\otimes^k}\right], \forall h\in H, \forall i\in\NN\right\}\\
&\cap& \left\{w\in \W:\sum_{k=1}^\infty
  \frac{x^k}{(k+1)!}\left|\nabla^{k+1} P_{s_i}
      f(w+h)\right|_{X\otimes H^{\otimes^{k+1}}}<\infty,\forall h\in
    H, \forall x\in\RR_+, \forall i\in\NN\right\}\eeaa

\no Observe that we know from \cite{BOOK} that $\left\{w\in \W:\sum_{k=1}^\infty \frac{x^n}{(k+1)!}\left|\nabla^{k+1}
    P_sf(w)\right|_{H^{\otimes^{k+1}}}<\infty, \forall
  x\in\RR_+\right\}$ and $\left\{w\in \W:P_{s_i}\nabla f(w+h)=P_{s_i}\nabla
  f(w)+\sum_{k=1}^\infty\frac{1}{k!}\nabla^k
  P_{s_i}\nabla f(w)\left[h^{\otimes^k}\right], \forall h\in H\right\}$
are of full measure  and H-invariant.\newline

\no Set $w\in A'$, $i\in\NN$, $h\in H$ and $h'\in H_0$ such that
$|h-h'|\leq\frac{1}{m}$ and $n\in\NN$ such that
$B\left(h,\frac{1}{m}\right)\subset B_n$,we have:

\beaa |P_{s_i}\nabla f(w+h)-\nabla f(w+h)|_{X\otimes H}&\leq&
|P_{s_i}\nabla f(w+h)-P_{s_i}\nabla f(w+h')|_{X\otimes H}\\
&&+ |P_{s_i}\nabla f(w+h')-\nabla f(w+h')|_{X\otimes H}\\
&&+ |\nabla f(w+h')-\nabla f(w+h)|_{X\otimes H}\\
&\leq&
|P_{s_i}((\nabla f(.+e^{-s_i})h)(w)-P_{s_i}(\nabla f(.+e^{-s_i}h'))(w)|_{X\otimes H}\\
&&+ |P_{s_i}\nabla f(w+h')-\nabla f(w+h')|_{X\otimes H}\\
&&+ |\nabla f(w+h')-\nabla f(w+h)|_{X\otimes H}\\
&\leq& P_{s_i}\theta_{nm}+ |P_{s_i}\nabla f(w+h')-\nabla f(w+h')|_{X\otimes H}+\theta_{nm}\eeaa

\no This proves that:
\beaa
&&\lim_{i\rightarrow\infty}
P_{s_i}\nabla f(w+h)=\nabla f(w+h)\;\;\forall h\in
H\eeaa
\no Now observe that for $w\in \W$ such that $(\theta_{n,m}(w))_{m\in\NN}$ converges
toward 0, $h\mapsto\nabla f(w+h)$ is uniformly continuous on $B_n$ hence
bounded. So for $h,k\in B_n$:
\beaa |f(w+h)-f(w+k)|_X&=&\left|\int_0^1\nabla f(w+\lambda
  h+(1-\lambda) k)[k-h]dt\right|\\
&\leq&\sup_{h'\in B_n}|\nabla f(w+h')|_{X\otimes H}|h-k|_H\eeaa
\no So the hypothesis imply that
$$\lim_{\epsilon\rightarrow 0}\sup_{h,k\in B_n,|h-k|_H\leq\epsilon}
\left| f(w+h)-f(w+k)\right|_{X}=0\;\; a.s.$$
where this supremum is a measurable random variable since $f\in \mathbb{D}_{p,1}(X)$
\no and we can construct a full measure $A''\subset \W$ similar to A'
where f takes the role of $\nabla f$. We denote $\tilde{A}=A'\cap
A''$. We have $\mu(\tilde{A})=1$ and for every $w\in\tilde{A}$ and
$h\in H$:
\beaa &&\lim_{i\rightarrow\infty}
P_{s_i} f(w+h)= f(w+h)\\
&&\lim_{i\rightarrow\infty}
P_{s_i}\nabla f(w+h)=\nabla f(w+h)\eeaa

\no Set $i,j\geq max(i_0,i_1,i_2)$, we have:

\no Now we can prove the differentiability of $h\mapsto f(w+h)$. Set $w\in \tilde{A}$ and $h\in H$, we aim to prove that:
$$\lim_{h'\rightarrow
  0}\frac{1}{|h'|_H}\left|f(w+h+h')-f(w+h)-\nabla
  f(w+h)[h']\right|_X=0$$

\no Set $h'\in H$, we have:

\beaa
&&\frac{1}{|h'|_H}\left|f(w+h+h')-f(w+h)-\nabla
  f(w+h)[h']\right|_X\\
&\leq&\frac{1}{|h'|_H}\left|f(w+h+h')
  -f(w+h)-\left(P_{s_i}f(w+h+h')
    -P_{s_i}f(w+h)\right)\right|_X\\
&&+\frac{1}{|h'|_H}\left|P_{s_i}f(w+h+h')
    -P_{s_i}f(w+h)-\nabla
      P_{s_i}f(w+h)[h']\right|_X\\
&&+\frac{1}{|h'|_H}\left|\nabla
      P_{s_i}f(w+h)[h']-\nabla f(w+h)[h']\right|_X\eeaa

\no We denote these three terms $A_{h'},B_{h'}$ and $C_{h'}$and we
deal with each one of them separately.
\beaa C_{h'}&\leq&\left|\nabla
      P_{s_i}f(w+h)-
      \nabla f(w+h)\right|_{X\otimes H}\\
&\leq&\left|e^{-s_i}
      P_{s_i}\nabla f(w+h)-
      \nabla f(w+h)\right|_{X\otimes H}\\
&\rightarrow& 0\eeaa

\no now the second term:

\beaa B_{h'}&=&\frac{1}{|h'|_H}\left|\sum_{k=2}^\infty
\frac{1}{k!}\nabla^k P_{s_i}f(w+h)\left[h'^{\otimes^k}
\right]\right|\\
&\leq&\frac{1}{|h'|_H}\sum_{k=2}^\infty
\frac{1}{k!}\left|\nabla^k P_{s_i}f(w+h)\right|_{X\otimes H^{\otimes^k}}\left|h'\right|_H^k\\
&=&\sum_{k=1}^\infty
\frac{1}{(k+1)!}\left|\nabla^{k+1} P_{s_i}f(w+h)\right|_{X\otimes H^{\otimes^{k+1}}}\left|h'\right|_H^k\\
&\rightarrow&0\eeaa

\no since $\sum_{k=1}^\infty
  \frac{1}{(k+1)!}\left|\nabla^{k+1} P_{s_i}f(w+h)\right|_{X\otimes
    H^{\otimes^{k+1}}}<\infty$.\newline
\beaa &&A_{h'}\\
&=&\lim_{j\rightarrow \infty}\frac{1}{|h'|_H}\left|P_{s_j}f(w+h+h')
  -P_{s_j}f(w+h)-\left(P_{s_i}f(w+h+h')
    -P_{s_i}f(w+h)\right)\right|_X\eeaa

\no We have
\beaa &&\left|P_{s_j}f(w+h+h')
  -P_{s_j}f(w+h)-\left(P_{s_i}f(w+h+h')
    -P_{s_i}f(w+h)\right)\right|_X\\
&\leq& \sup_{\lambda\in [0,1]}\left|\nabla
  P_{s_j}f(w+h+\lambda h')
  -\nabla P_{s_i}f(w+\lambda
  h+\lambda h')\right|_{X\otimes H}|h'|_H\eeaa

\no Now set $\lambda\in [0,1]$, $\epsilon>0$. We can assume that $|h-h'|_H\leq\frac{1}{m}$, set $n\in\NN$ such that
$B\left(h,\frac{1}{m}\right)\subset B_n$.
We have:
\beaa &&\left|\nabla
  P_{s_j}f(w+h+\lambda h')
  -\nabla P_{s_i}f(w+
  h+\lambda h')\right|_{X\otimes H}\\
&\leq& \left|\nabla
  P_{s_j}f(w+h+\lambda h')
  -\nabla P_{s_j}f(w+h)\right|_{X\otimes H}\\
&&+\left|\nabla P_{s_j}f(w+h)-\nabla P_{s_i}f(w+h)\right|_{X\otimes H}\\
&&+\left|\nabla
  P_{s_i}f(w+h)
  -\nabla P_{s_i}f(w+\lambda
  h+(1-\lambda)h')\right|_{X\otimes H}\\
&\leq& \left|e^{-s_j}
  P_{s_j}(\nabla f(.+e^{-s_j}(h+\lambda h')))(w)
  -e^{-s_j}P_{s_j}(\nabla f(.+e^{-s_j}h))(w)\right|_{X\otimes H}\\
&&+\left|\nabla P_{s_j}f(w+h)-\nabla P_{s_i}f(w+h)\right|_{X\otimes H}\\
&&+\left|e^{-s_i}
  P_{s_i} (\nabla f(.+e^{-s_i}h))(w)
  -e^{-s_i}P_{s_i}(\nabla f(.+e^{-s_i}(
  h+\lambda h')))(w)\right|_{X\otimes H}\\
&\leq&P_{s_i}\theta_{nm}(w+h)+P_{s_j}\theta_{nm}(w+h)+\left|\nabla
  P_{s_j}f(w+h)-\nabla P_{s_i}f(w+h)\right|_{X\otimes H}\eeaa
\no which is smaller than $\epsilon$ for i and j large enough. It
ensures that $A_{h'}$ tends toward 0 when h' converges toward 0, which
concludes the proof.\nqed

\section{\bf{Extension to higher order derivatives}}

\begin{corollary}
\no Assume that $f:\W\rightarrow X$ is in $\mathbb{D}_{p,r}(X)$ for
some $p>1$. Assume that
$h\mapsto \nabla^k f(w+h)$ is $\mu$-a.s. uniformly continuous on every $B_n$.\newline
\no Then f is $H-C^r$ and its H-derivatives up to order n are equal to
its weak derivatives of the same order.
\end{corollary}

\nproof We prove this with a recurrence over n. The case $r=1$ is
theorem \ref{th}. Now set $r\leq 2$ and assume that the result is proven for
every integer up to k-1. Set $n\in\NN$ and A a measurable subset of $\W$ such that $\mu(A)=1$ and for every $w\in A$ $h\mapsto \nabla^r f(w+h)$ is uniformly continuous on $B_n$. Set $w\in A$, $B_n$ being closed, $h\mapsto \nabla^r f(w+h)$ is bounded is bounded on $B_n$. Consequently,  $h\mapsto\nabla^{r-1}f(w+h)$ is lipschitz on $B_n$ and so is uniformly continuous on $B_n$. The recurrence hypothesis ensures that f is H-$C^{r-1}$ and that its H derivatives up to order r-1 are equals to its weak derivatives of the same order.\newline
\no Applying theorem \ref{th} to $\nabla^{r-1} f$, we get that it is H-$C^1$ and that its H-derivative is $\nabla^r f$, which conclude the proof.\nqed

\vspace{2cm}
\footnotesize{
\noindent
K\'evin HARTMANN, Institut Telecom, Telecom ParisTech, LTCI CNRS D\'ept. Infres, \\
23 avenue d'Italie, 75013, Paris, France\\
kevin.hartmann@polytechnique.org}

\end{document}